\documentclass[reqno]{amsart}
\newtheorem{thm}{Theorem}[section]

\newtheorem{prop}{Proposition}[section]

\def \rn {\mathbb{R}^4}
\def \nn {\mathbb{N}}
\def \rr {\mathbb{R}}

\def \tuk {\tilde{u}_k}
\def \rk {r_k}
\def \sk {s_k}
\def \xk {x_k}
\def \yk {y_k}
\def \wk {w_k}
\def \Uk {U_k}
\def \xjk {x_{j,k}}
\def \xik {x_{i,k}}
\def \ovb {\overline{B}}
\def \rhok {\rho_{k}}
\def \txj {\tilde{x}_{j}}
\def \txjk {\tilde{x}_{j,k}}
\def \txi {\tilde{x}_{i}}
\def \txik {\tilde{x}_{i,k}}
\def \mik {\mu_{i,k}}
\def \Vk {V_k}
\def \tVk {\tilde{V}_k}
\def \uk {u_k}
\def \mk {\mu_k}
\def \vk {v_k}
\def \Rk {R_k}

\def \txpik {\tilde{x}_{\phi(i),k}}
\def \txpi {\tilde{x}_{\phi(i)}}
\def \txpj {\tilde{x}_{\phi(j)}}
\def \xpik {x_{\phi(i),k}}
\def \beq {\begin{eqnarray*}}
\def \eeq {\end{eqnarray*}}
\def \beqn {\begin{eqnarray}}
\def \eeqn {\end{eqnarray}}
\def \bequa {\begin{equation}}
\def \eequa {\end{equation}}
\title[Quantization]{Quantization effects for a fourth order equation of exponential growth in dimension four}
\author{Fr\'ed\'eric Robert}
\email{frobert@math.unice.fr}
\address{Universit\'e de Nice-Sophia Antipolis, Laboratoire J.A.Dieudonn\'e, Parc Valrose, 06108 Nice Cedex 2, France}
\date{December 7th, 2005.}
\begin{document}

\begin{abstract}
We investigate 
the asymptotic behavior as $k \to +\infty$ of 
sequences $(\uk)_{k\in\nn}\in C^4(\Omega)$ of solutions of 
the equations $\Delta^2\uk=\Vk e^{4\uk}$ on $\Omega$, where $\Omega$ is a bounded domain of $\rn$ and $\lim_{k\to +\infty}\Vk=1$ in
$C^0_{loc}(\Omega)$. The corresponding 2-dimensional problem was
studied by Br\'ezis-Merle and Li-Shafrir who pointed out 
that there is a quantization of the energy when blow-up occurs. As shown by Adimurthi, Struwe and the
author
\cite{adirost}, such a quantization does not hold in dimension four for the problem 
in its full generality. We prove here that under
natural hypothesis on $\Delta\uk$, we recover such a quantization as in dimension 2.
\end{abstract}

\maketitle
\section{Introduction}
Let $\Omega$ be a bounded domain of $\rn$. Let a sequence $(\Vk)_{k\in\nn}\in C^0(\Omega)$ such that
\begin{equation}\label{lim:Vk}
\lim_{k\to +\infty}\Vk=1
\end{equation}
in $C^0_{loc}(\Omega)$. Let $(\uk)_{k\in\nn}$ be a sequence of functions in $C^4(\Omega)$ such that
$$\Delta^2\uk= \Vk e^{4\uk}\eqno{(E)}$$
in $\Omega$ for all $k\in\nn$. Here and in the sequel, $\Delta=-\sum\partial_{ii}$ is the Laplacian with minus sign convention. 
In this paper, we address the question of the asymptotics of the $\uk$'s when $k\to +\infty$. A natural (and simple) behavior is when 
there exists
$u\in C^4(\Omega)$ such that 
\begin{equation}\label{eq:cv}
\lim_{k\to +\infty}\uk=u
\end{equation}
in $C^{3}_{loc}(\Omega)$. In this situation, we say that $(\uk)_{k\in\nn}$ is relatively compact in $C^3_{loc}(\Omega)$. 
However, the structure of equation $(E)$ is much richer due to its scaling invariance properties. 
The scaling invariance is as follows. Given $k\in\nn$, $x_k\in\Omega$
and $\mu_k>0$ , we let
\begin{equation}\label{eq:inv:uk}
\tilde{u}_{k}(x):=\uk(\xk+\mu_k x)+\ln\mu_k
\end{equation}
for all $x\in\mu_k^{-1}(\Omega-\xk)$. Letting $\tilde{V}_k(x)=\Vk(x_k+\mu_k x)$ for all $x\in \mu_k^{-1}(\Omega-\xk)$, we get 
that the rescaled function $\tilde{u}_{k}$ satisfies
$$\Delta^2\tilde{u}_{k}=\tilde{V}_k e^{4\tilde{u}_{k}}$$
on $\mu_k^{-1}(\Omega-\xk)$ -- an equation like $(E)$. This scaling invariance forces some situations more subtle that (\ref{eq:cv}) to
happen. A very simple example is the following: we consider a sequence $(\mu_k)_{k\in\nn}\in\rr_+^\star$ such that $\lim_{k\to +\infty}\mu_k=0$
and for any $k\in\nn$, we define the function
$$f_k(x)=\ln \frac{\sqrt{96}\mk}{\sqrt{96}\mk^2+|x|^2}$$
for all $x\in\rn$. Then $f_k$ satisfies $(E)$ with $\Vk\equiv 1$ for all $k\in\nn$. The sequence 
$(f_k)_{k\in\nn}$ does not converge in $C^{0}_{loc}(\rn)$: we have that
$$\lim_{k\to +\infty}f_k(0)=+\infty\hbox{ and }\lim_{k\to +\infty}f_k=-\infty\hbox{ uniformly locally on }\rn\setminus\{0\}.$$
In addition, we get that
$$\Vk e^{4 f_k}\, dx\rightharpoonup 16\pi^2\delta_0$$
when $k\to +\infty$ weakly for the convergence of measures. Scaling as in (\ref{eq:inv:uk}), we get that
$$\lim_{k\to +\infty}f_k(\mu_k x)-f_k(0)=\ln\frac{\sqrt{96}}{\sqrt{96}+|x|^2}$$
for all $x\in\rn$. Concerning terminology, we say that the sequence $(\uk)_{k\in\nn}$ blows-up if it is not relatively 
compact in $C^3_{loc}(\Omega)$, so that, up to any subsequence, (\ref{eq:cv}) does not hold. In the above 
example, the $(f_k)$'s blow up. In this paper,
we are mainly concerned with the blow-up behavior of solutions of $(E)$.

\medskip\noindent In dimension two, the corresponding problem has been
studied (among others) by Br\'ezis-Merle \cite{brezismerle} and
Li-Shafrir \cite{lishafrir}. We also refer to Druet \cite{druet} and
Adimurthi-Struwe \cite{adistruwe} for the description of equations with more intricate nonlinearities and to Tarantello \cite{t} for equations with singularities. Li and Shafrir proved the following:

\begin{thm}[Li-Shafrir \cite{lishafrir}]\label{th:ls} Let $\Sigma$ be a bounded domain of $\rr^2$, 
$(\bar{V}_k)_{k\in\nn}\in C^0(\Sigma)$ be a
sequence of functions such that $\lim_{k\to +\infty}\bar{V}_k=1$ in $C^0_{loc}(\Sigma)$, and 
$(\bar{u}_k)_{k\in\nn}\in C^2(\Sigma)$ be a sequence 
such that 
$$\Delta\bar{u}_k=\bar{V}_k e^{2\bar{u}_k}$$
in $\Sigma$ for all $k\in\nn$, and such that there exists $\Lambda\in\rr$ such that 
$\int_{\Sigma}\bar{V}_k e^{2\bar{u}_k}\, dx\leq \Lambda$ 
for all $k\in\nn$. Then either (i) the sequence $(\uk)_{k\in\nn}$ is relatively compact in $C^1(\Omega)$, or (ii) there exists $N\in\nn$, there exist $\bar{x}_1,...,\bar{x}_N\in\Omega$, there exist $\bar{\alpha}_1,...,\bar{\alpha}_N\in \nn^\star$ such that, up to a subsequence
$$\bar{V}_ke^{2\bar{u}_k}\rightharpoonup \sum_{i=1}^N4\pi \bar{\alpha}_i\delta_{\bar{x}_i}$$
weakly for the convergence of measures when $k\to +\infty$. Moreover, $\lim_{k\to +\infty}\bar{u}_k=-\infty$ uniformly locally in $\Sigma\setminus\{\bar{x}_1,...,\bar{x}_N\}$.
\end{thm}
We refer to this statement as a quantization result. The justification of this terminology is as follows: if in Theorem \ref{th:ls} we have blow-up (that is case (i) does not hold), then for any $\omega\subset\subset\Sigma$ such that $\partial\omega\cap\{\bar{x}_1,...,\bar{x}_N\}=\emptyset$, we have that $\lim_{k\to +\infty}\int_\omega\bar{V}_k e^{2\bar{u}_k}\, dx\in 4\pi\nn$. Moreover, the sequence $(\uk)_{k\in\nn}$ develop singularities on a set at most finite, that is $\{\bar{x}_1,...,\bar{x}_N\}$.

\medskip\noindent Surprisingly, such a quantization result is false when we come back to our initial four-dimensional problem $(E)$. In a joint work with Adimurthi and Michael Struwe \cite{adirost}, we exhibit a sequence of solutions to $(E)$ that blows-up, carry a non-quantified energy and develop singularities on a hypersurface of $\rn$. In \cite{adirost}, we described the behaviour of arbitrary solutions to $(E)$ and proved that any blowing-up sequence $(\uk)_{k\in\nn}$ concentrates at the zero set of a nonpositive bi-harmonic function, and that outside this set, $\lim_{k\to +\infty}\uk=-\infty$ uniformly. In view of the examples provided in \cite{adirost}, this result is optimal. Therefore, giving a more precise description requires additional hypothesis on $(\uk)$. 

\medskip\noindent A natural hypothesis is to impose a Navier boundary condition, (that is $\uk=\Delta\uk=0$ on $\partial\Omega$) or a Dirichlet condition (that is $\uk=\frac{\partial\uk}{\partial\nu}=0$ on $\partial\Omega$): actually, in these cases, we get that there is no blow-up and we recover relative compactness. Wei \cite{wei} studied a problem similar to $(E)$ assuming that $\Delta\uk=0$ on $\partial\Omega$ and $\uk=c_k$ on $\partial\Omega$, where $(c_k)_{k\in\nn}\in\rr$ is a sequence of real numbers such that $\lim_{k\to +\infty}c_k=-\infty$: in this context, Wei describes precisely the asymptotics and recovers quantization. Another natural hypothesis is to assume that the functions $\uk$ are radially symmetrical: in this situation, we describe completely the asymptotics in \cite{rradial}. In all these situations, the critical quantity to observe happens to be $\Delta\uk$ as shown in the following example. We let $\alpha\in (0,16\pi^2)$. It follows from \cite{changchen} that there exists $v\in C^4(\rn)$ radially symmetrical such that $v\leq v(0)=0$ and
$$\Delta^2 v= e^{4v}\hbox{ in }\rn\hbox{ and }\int_{\rn}  e^{4v}\, dx=\alpha.$$
Contrary to the two-dimensional case, where the only solutions to the corresponding equation are of a type similar to $f_k$ with a quantization of the energy, we get in four dimensions many solutions with arbitrary small energy. More precisely, it follows from \cite{lin} that there exists $C>0$ such that $\Delta v(x)\geq C$ for all $x\in\rn$. For any $k\in\nn^\star$, we define the function
$$g_k(x)=v(kx)+\ln k$$
for all $x\in\rn$. As easily checked, due to the scaling invariance (\ref{eq:inv:uk}) of $(E)$, $g_k$ verifies $(E)$ with $\Vk\equiv 1$. We also get that the sequence $(g_k)_{k\in\nn}$ blows up. It follows from straightforward computations that
$$\lim_{k\to +\infty}\int_{B_{1}(0)}\Vk e^{4g_k}\, dx=\alpha.$$
Moreover, for any $\omega\subset\rn$, we have that 
$$\lim_{k\to +\infty}\Delta g_k=+\infty$$
uniformly in $\omega$. Since $\alpha>0$ can be chosen as small as we want, we then get blowing-up sequences with arbitrary positive small energy, and there is no quantization here. Note that concerning the sequence $(f_k)_{k\in\nn}$ of the first example, we have that for any $\omega\subset\subset\rn\setminus\{0\}$, there exists $C(\omega)>0$ such that 
$$|\Delta f_k(x)|\leq C(\omega)$$
for all $k\in\nn^\star$ and all $x\in \omega$. The fundamental difference between the $(f_k)$'s and the $(g_k)$'s is that in the first case, the Laplacian is bounded outside the singularity, and in the second case, the Laplacian goes to $+\infty$ uniformly. This fact is actually general. 
The objective of this paper is to prove the following result:

\begin{thm}\label{th:intro} Let $\Omega$ be a bounded domain of $\rn$, 
$(\Vk)_{k\in\nn}\in C^0(\Omega)$ be a sequence such that (\ref{lim:Vk}) holds, and 
$(\uk)_{k\in\nn}$ be a sequence of functions in $C^4(\Omega)$ such that $(E)$ holds, and such that there exists $\Lambda>0$ such that
$\int_\Omega \Vk e^{4\uk}\, dx\leq \Lambda$ 
for all $k\in\nn$. Assume there exist $C>0$ and $\omega_0\subset\subset\Omega$ such that
$\Vert(\Delta\uk)_-\Vert_1\leq C$ and 
$$\Vert \Delta\uk\Vert_{L^1(\omega_0)}\leq C$$ 
for $k\in\nn$. Then (i) either $(\uk)_{k\in\nn}$ is relatively compact in $C^3_{loc}(\Omega)$, or (ii) there exists $N\in\nn$, there exist $x_1,...,x_{N}\in\Omega$, there exist $\alpha_1,...,\alpha_{N}\in\mathbb{N}^\star$ such that
$$\Vk e^{4\uk}\rightharpoonup \sum_{i=1}^{N}16\pi^2\alpha_i\delta_{x_i}$$
weakly in the sense of measures when $k\to +\infty$ up to a subsequence. Moreover, still in Case (ii), we have that
$\lim_{k\to +\infty}\uk=-\infty$ uniformly locally in $\Omega\setminus\{x_1,...,x_{N}\}$.
\end{thm}

As a remark, note that the control of the positive part of $\Delta\uk$
is only required on an arbitrary subdomain of $\Omega$. This result is
optimal as shown in the preceding example involving the function
$g_k$. In a joint work with Olivier Druet \cite{druetrobert}, we
studied the corresponding problem on four-dimensional Riemannian
manifolds, where the bi-Laplacian is replaced by a fourth-order
elliptic operator refered to as $P$: when the kernel of $P$ is such
that $Ker\, P=\{constants\}$, we get similar results as in Theorem
\ref{th:intro} with the additional information that $\alpha_i=1$ for
all $i\in \{1,...,N\}$. The techniques used in \cite{druetrobert} are
different from the techniques used here: the main reason is that for
equation $(E)$, the kernel of the bi-Laplacian contains more than the
constant functions. Related references in the context of Riemannian
manifolds are Malchiodi \cite{malchiodi} and Malchiodi-Struwe \cite{malchiodistruwe}. As a
remark, the corresponding question in dimension $n\geq 5$ was
considered in Hebey-Robert \cite{hebeyrobert}.

\medskip\noindent This paper is organized as follows. In section
\ref{sec:construc}, we prove that under our hypothesis, concentration
holds on finitely many points and not on a hypersurface. In section
\ref{sec:est}, we prove that, up to rescaling, the $\uk$'s converge to
a generic pattern when $k\to +\infty$. In section \ref{sec:blowup}, we
analyse precisely the blow-up and we prove Theorem \ref{th:intro} in section \ref{sec:conc}. In the sequel, 
$C$ denotes a positive constant, with value allowed to change from one line to the other. Note also that all the 
convergence results are up to a
subsequence, even when it is not precised.

\medskip\noindent{\bf\large Acknowledgement:} the author thanks Adimurthi 
and Michael Struwe for having pointed out this problem, and also thanks them for stimulating discussions. The author thanks Emmanuel Hebey for stimulating discussions on this problem.

\section{Construction of the concentration points}\label{sec:construc}
In the sequel, we let $\Omega$ be a bounded domain of $\rn$. We let a sequence a sequence $(\Vk)_{k\in\nn}\in C^0(\Omega)$ such that (\ref{lim:Vk}) holds. Let $(\uk)_{k\in\nn}$ be a sequence of functions in $C^4(\Omega)$ such that $(E)$ holds. We assume that there exists $\Lambda>0$ such that 
\bequa\label{hyp:lambda}
\int_\Omega \Vk e^{4\uk}\, dx\leq \Lambda
\eequa
for all $k\in\nn$. We assume that there exist $\omega_0\subset\subset\Omega$ and $C>0$ such that
\bequa\label{hyp:ellip}
\Vert \Delta\uk\Vert_{L^1(\omega_0)}\leq C
\eequa
and
\bequa\label{hyp:min}
\Vert(\Delta\uk)_-\Vert_1\leq C
\eequa
for all $k\in\nn$. The objective of this section is to prove that the $(\uk)$'s concentrate at a finite number of points. This is the object of the following proposition:

\begin{prop}\label{bound:out0} Let $\Omega$ be a bounded domain of $\rn$. Let a sequence $(\Vk)_{k\in\nn}\in C^0(\Omega)$ such that (\ref{lim:Vk}) holds. Let $(\uk)_{k\in\nn}$ be a sequence of functions in $C^4(\Omega)$ such that $(E)$ holds. We assume that there exists $\Lambda>0$ such that \eqref{hyp:lambda} holds. We assume that \eqref{hyp:min} and \eqref{hyp:ellip} hold. We let 
\bequa\label{def:S0}
S_0:=\left\{x\in\Omega/ \liminf_{\delta\to 0}\liminf_{k\to +\infty}\int_{B_\delta(x)}\Vk  e^{4\uk}\, dy\geq 8\pi^2\right\}.
\eequa
which is a finite set. Then for any $\omega\subset\subset \Omega\setminus S_0$, there exists $C(\omega)>0$ such that
$$|\Delta\uk(x)|\leq C(\omega)\hbox{ and }\uk(x)\leq C(\omega)$$
for all $x\in\omega$ and all $k\in\nn$. More precisely, we are in one
and only one of the following situations:

\smallskip\noindent (A1) there exists $u\in C^4(\Omega\setminus S_0)$ such that 
$$\lim_{k\to +\infty}\uk=u\hbox{ in }C^3_{loc}(B_\delta(x_0))$$

\noindent (A2) $\lim_{k\to +\infty}\uk=-\infty$ uniformly locally on
$\Omega\setminus S_0$.

\end{prop}

The proof of Proposition \ref{bound:out0} proceeds in two steps. Note that it follows from \eqref{hyp:lambda} that $S_0$ is at most finite. We let a sequence $(\Vk)_{k\in\nn}\in C^0(\Omega)$ such that (\ref{lim:Vk}) holds. Let $(\uk)_{k\in\nn}$ be a sequence of functions in $C^4(\Omega)$ such that $(E)$ holds. We assume that there exists $\Lambda>0$ such that \eqref{hyp:lambda} holds. We assume that \eqref{hyp:min} and \eqref{hyp:ellip} hold.

\medskip\noindent{\bf Step \ref{sec:construc}.1:} We let $x_0\in\Omega\setminus S_0$. We claim that we are in one and only one of the following situations:

\smallskip\noindent (B1) there exists $u\in C^4(B_\delta(x_0))$ such that $\lim_{k\to +\infty}\uk=u$ in $C^3_{loc}(B_\delta(x_0))$.\par
\noindent (B2) there exists $\phi\in C^4(B_\delta(x_0))$ such that $\Delta^2\phi=0$, $\phi\leq 0$, $\phi\not\equiv 0$ there exists a sequence $(\beta_k)_{k\in\nn}\in\rr_+^\star$ such that $\lim_{k\to +\infty}\beta_k=+\infty$ and 
$$\lim_{k\to +\infty}\frac{\uk}{\beta_k}=\phi$$
in $C^3_{loc}(B_\delta(x_0))\cap\{\phi<0\}$. 

\medskip\noindent{\it Proof of the claim:} This claim is a particular case of the Theorem obtained in \cite{adirost}. As a preliminary remark, note that the two cases (B1) and (B2) are disjoint. Since $x_0\in\Omega\setminus S_0$, we let $\delta>0$ and $\alpha<8\pi^2$ such that
$$\int_{B_\delta(x_0)}\Vk  e^{4\uk}\, dx\leq \alpha<8\pi^2$$
for all $k\in\nn$. We let $\vk$ such that
\bequa\label{def:vk}
\Delta^2\wk=\Vk  e^{4\uk}\hbox{ in }B_\delta(x_0),\quad \wk=\Delta\wk=0\hbox{ on }\partial B_\delta(x_0).
\eequa
It follows from \cite{lin} (see also \cite{wei}) that there exists $p>1$ such that
\bequa\label{ineq:nrjve}
\int_{B_\delta(x_0)} e^{4p|\wk|}\, dx\leq C
\eequa
for all $k\in\nn$. We let $h_k:=\uk-\wk$ on $B_\delta(x_0)$. Clearly $\Delta^2 h_k=0$. It follows from \eqref{hyp:lambda} and \eqref{ineq:nrjve} that $\Vert (h_k)_+\Vert_{L^1(B_\delta(x_0))}=O(1)$ when $k\to +\infty$. We distinguish two situations:

\smallskip\noindent{\it Case \ref{sec:construc}.1.1:} We assume that $\Vert h_k\Vert_{L^1(B_{\delta/2}(x_0))}=O(1)$ when $k\to +\infty$. Since $h_k$ is bi-harmonic, there exists $h_\infty\in C^4(B_\delta(x_0))$ such that 
\bequa\label{cv:case:1}
\lim_{k\to +\infty}h_k=h_\infty
\eequa
in $C^4_{loc}(B_\delta(x_0))$. We refer to \cite{adirost} for details about this assertion. Plugging \eqref{ineq:nrjve} and \eqref{cv:case:1} in \eqref{def:vk}, we get that $(\wk)_{k\in\nn}$ is bounded in $C^0_{loc}(B_\delta(x_0))$, and so is $(\uk)_{k\in\nn}$. It then follows from standard elliptic theory that there exists $u\in C^4(B_\delta(x_0))$ such that $\lim_{k\to +\infty}\uk=u$ in $C^3_{loc}(B_\delta(x_0))$, and we recover Case (B1) of the claim. This proves the claim in Case \ref{sec:construc}.1.1.

\smallskip\noindent{\it Case \ref{sec:construc}.1.2:} We assume that $\lim_{k\to +\infty}\Vert h_k\Vert_{L^1(B_{\delta/2}(x_0))}=+\infty$. Since $h_k$ is bi-harmonic, there exists $\phi\in C^4(B_\delta(x_0))\setminus\{0\}$ such that $\Delta^2\phi=0$, $\phi\leq 0$, there exists a sequence $(\beta_k)_{k\in\nn}\in\rr_+^\star$ such that $\lim_{k\to +\infty}\beta_k=+\infty$ and such that 
\bequa\label{cv:case:2}
\lim_{k\to +\infty}\frac{h_k}{\beta_k}=\phi
\eequa
in $C^4_{loc}(B_\delta(x_0))$. We refer to \cite{adirost} for details about this assertion. In particular, $h_k\to -\infty$ uniformly locally on $\phi<0$. Arguing as in Case \ref{sec:construc}.1.1, we then obtain that $(\wk)_{k\in\nn}$ converges in $C^3_{loc}(B_{\delta}(x_0)\cap\{\phi<0\})$. It then follows from \eqref{cv:case:2} that $\lim_{k\to +\infty}\frac{\uk}{\beta_k}=\phi$ in $C^3_{loc}(B_{\delta}(x_0)\cap\{\phi<0\})$, and we recover Case (B2) of the claim. This proves the claim in Case \ref{sec:construc}.1.2.\hfill$\Box$

\medskip\noindent{\bf Step \ref{sec:construc}.2:} We are in position to prove Proposition \ref{bound:out0}. Since $\Omega\setminus S_0$ is connected and harmonic functions are analytic, it follows from Step \ref{sec:construc}.1 that we are in one and only one of the following situations:

\medskip\noindent{\it Case \ref{sec:construc}.2.1:} There exists $u\in C^4(\Omega\setminus S_0)$ such that $\lim_{k\to +\infty}\uk=u$ in $C^3_{loc}(\Omega\setminus S_0)$. In this situation, we recover Case (A1) of Proposition \ref{bound:out0}.

\medskip\noindent{\it Case \ref{sec:construc}.2.2:} There exists $\phi\in C^4(\Omega\setminus S_0)$ such that $\Delta^2\phi=0$, $\phi\leq0$, $\phi\not\equiv 0$, there exists a sequence $(\beta_k)_{k\in\nn}\in\rr_+^\star$ such that $\lim_{k\to +\infty}\beta_k=+\infty$ and 
\bequa\label{lim:step2:2}
\lim_{k\to +\infty}\frac{\uk}{\beta_k}=\phi\hbox{ in }C^3_{loc}(\Omega\cap\{\phi<0\}\setminus S_0).
\eequa
We claim that $\Delta\phi\equiv 0$. Indeed, there exists $x\in\omega_0$ ($\omega_0$ was defined in \eqref{hyp:ellip}) such that $\phi(x)<0$ (otherwise $\phi\equiv 0$ on $\omega_0$ and then $\phi\equiv 0$ on $\Omega\setminus S_0$ since harmonic fonctions are analytic. A contradiction). We then get that \eqref{lim:step2:2} holds in a neighborhood of $x_0$. With \eqref{hyp:ellip}, we then get that $\Delta\phi=0$ in a neighborhood of $x$. Since $\Delta\phi$ is harmonic, and therefore analytic, we get that $\Delta\phi\equiv 0$ on $\Omega\setminus S_0$. This proves the claim.

\medskip\noindent Since $\phi\not\equiv 0$, $\phi\leq 0$ and $\Delta\phi=0$, it follows from the maximum principle that $\phi<0$ on $\Omega\setminus S_0$. Consequently, 
$$\lim_{k\to +\infty}\frac{\uk}{\beta_k}=\phi\hbox{ in }C^3_{loc}(\Omega\setminus S_0).$$
In particular, we get that $\lim_{k\to +\infty}\uk=-\infty$ uniformly locally on $\Omega\setminus S_0$. With the equation $(E)$ and \eqref{hyp:min}, it follows from elliptic theory that either $\lim_{k\to +\infty}\Delta\uk=+\infty$ uniformly locally in $\Omega\setminus S_0$, or $(\Delta\uk)_{k\in\nn}$ is uniformly bounded when $k\to +\infty$ locally in $\Omega\setminus S_0$: it follows from hypothesis \eqref{hyp:ellip} that the first situation cannot hold, and we get that Case (B2) of Proposition \ref{bound:out0} holds.

\medskip\noindent Clearly Proposition \ref{bound:out0} is a consequence of Steps \ref{sec:construc}.1 and \ref{sec:construc}.2.

\section{Pointwise estimates}\label{sec:est}
This section is devoted to the proof of the following Proposition:
\begin{prop}\label{estimate:pointwise1} Let $\Omega$ be a bounded domain of $\rn$. We let a sequence $(\Vk)_{k\in\nn}\in C^0(\Omega)$ such that (\ref{lim:Vk}) holds. Let $(\uk)_{k\in\nn}$ be a sequence of functions in $C^4(\Omega)$ such that $(E)$ holds. We assume that there exists $\Lambda>0$ such that \eqref{hyp:lambda} holds. We assume that \eqref{hyp:min} and \eqref{hyp:ellip} hold. We assume that
$$S_0\neq\emptyset.$$
Then there exists $N\in\nn^\star$, there exists families of points $(x_{1,k})_{k\in\nn},..., (x_{N,k})_{k\in\nn}$ in $\Omega$ such that for all $i\in \{1,...,N\}$, we have that $\lim_{k\to +\infty}x_{i,k}=x_i\in S_0$ and such that for any $\omega\subset\subset\Omega$, there exists $C(\omega)>0$ such that
$$(\inf_{i\in\{1,...,N\}} |x-\xik|) e^{\uk(x)}\leq C(\omega)\hbox{ and }(\inf_{i\in\{1,...,N\}} |x-\xik|)^2|\Delta\uk(x)|\leq C(\omega)$$
for all $k\in\nn$ and all $x\in\omega$. Moreover,
$$\lim_{k\to +\infty}\frac{|x_{i,k} -x_{j,k}|}{ e^{-\uk(\xik)}} = + \infty \hbox{ for all } i\neq j,$$
and for any $i\in \{1,...,N\}$ and any $x\in\rn$, we have that
$$\lim_{k\to +\infty}(\uk(\xik+ e^{-\uk(\xik)}x)-\uk(\xik))=\ln\frac{\sqrt{96}}{\sqrt{96}+|x|^2}.$$
Moreover, this convergence holds in $C^3_{loc}(\rn)$. In addition, $\lim_{k\to +\infty}\uk=-\infty$ uniformly on every compact subset of $\Omega\setminus S_0$.
\end{prop}

This section is devoted to the proof of Proposition \ref{estimate:pointwise1}. We let $\omega\subset\subset\Omega$. Up to taking $\omega$ larger, we assume that $S_0\subset\omega$. We follow the proof of \cite{rs}. We let $\xk\in\overline{\omega}$ such that
$$\uk(\xk)=\sup_{\omega}\uk.$$
Since $S_0\neq\emptyset$ and $S_0\subset\omega$, we get that $\lim_{k\to +\infty}\uk(\xk)=+\infty$. In this situation, it follows from Proposition \ref{bound:out0} that $\lim_{k\to +\infty}\xk=x_0\in S_0\cap\overline{\omega}=S_0\cap\omega$. We let $\delta>0$ small such that
$$B_{2\delta}(\xk)\subset\omega$$
for all $k\in\nn$. We define 
\bequa\label{def:vk:proof}
\mk:=e^{-\uk(\xk)}\hbox{ and }\vk(x):=\uk(\xk+\mk x)-\uk(\xk)
\eequa
and $\tVk(x):=\Vk(\xk+\mk x)$ for $|x|<\frac{2\delta}{\mk}$ and all $k\in\nn$. Equation $(E)$ yields
\bequa\label{eq:vk:new}
\Delta^2\vk=\tVk e^{4\vk},
\eequa
with $\vk(x)\leq\vk(0)=0$. 

\medskip\noindent{\bf Step \ref{sec:est}.1:} We claim that there existe $C>0$ independant of $k$ and $R$ such that
\bequa\label{estim:intve}
\int_{B_R(0)}|\Delta\vk|\, dx\leq C R^2 +C R^4\mk^2
\eequa
for all $k\in\nn$ and all $R<\delta\mk^{-1}$. We prove the claim. We let $G_{\delta,k}$ be the Green's function for the Laplacian on $B_\delta(\xk)$ with Dirichlet boundary condition. We get that
$$\Delta\uk(z)=\int_{B_\delta(\xk)}G_{\delta,k}(z,y)\Delta^2\uk(y)\, dy +\varphi_k(z)$$
for all $z\in B_\delta(\xk)$, where $\varphi_k$ is the unique harmonic function on $B_\delta(\xk)$ such that $\varphi_k(y)=\Delta\uk(y)$ for all $y\in \partial B_{\delta}(\xk)$. With Proposition \ref{bound:out0} and the comparison principle, we get that there exists $C>0$ such that 
\bequa\label{ineq:phi:k}
|\varphi_k(z)|\leq C
\eequa 
for all $z\in B_\delta(\xk)$. We let $x\in \rn$ such that $|x|<\delta\mk^{-1}$. Using the definition \eqref{def:vk:proof} of $\vk$, we get that
$$\Delta\vk(x)= \int_{B_\delta(\xk)}\mk^2 G_{\delta,k}(\mk x,y)\Delta^2\uk(y)\, dy+\mk^2\varphi_k(\xk+\mk x).$$
Integrating this equation, using $(E)$, (\ref{hyp:lambda}), \eqref{ineq:phi:k} and standard estimates on the Green's function, we get that
\beq
\int_{B_R(0)}|\Delta\vk|\, dx& \leq & C \int_{x\in B_R(0)}\int_{y\in B_\delta(\xk)}\mk^2 G_{\delta,k}(\mk x,y)e^{4\uk(y)}\, dy\, dx+ C R^4 \mk^2\\
&\leq & C\int_{B_\delta(\xk)} e^{4\uk(y)}\left(\int_{B_R(0)}\frac{\mk^2}{|\mk x-y|^2}\, dx\right)\, dy+ C R^4 \mk^2\\
&\leq & C\int_{B_\delta(\xk)}e^{4\uk(y)}\left( CR^2\right)\, dy\leq C \Lambda R^2+ C R^4 \mk^2,
\eeq
for all $k\in\nn$. This proves the claim.

\medskip\noindent{\bf Step \ref{sec:est}.2:} We claim that for any $x\in\rn$, we have that 
\bequa\label{lim:vk:ellip}
\lim_{k\to +\infty}\vk(x)= \ln\frac{\sqrt{96}}{\sqrt{96}+|x|^2}:=U_0(x),
\eequa
moreover this convergence holds in $C^3_{loc}(\rn)$. We briefly prove the claim. With \eqref{estim:intve}, we get that $\Delta\vk$ is bounded in $L^1_{loc}$. Since $\vk\leq\vk(0)=0$, it then follows from \eqref{eq:vk:new} and standard elliptic theory that, up to a subsequence, there exists $v\in C^4(\rn)$ such that $\lim_{k\to +\infty}\vk=v$ in $C^3_{loc}(\rn)$, with $\Delta^2 v=e^{4v}$ and $e^{4v}\in L^1(\rn)$. Passing to the limit $k\to +\infty$ in \eqref{estim:intve} and using the classification of Lin \cite{lin}, we get that $v\equiv U_0$. We refer to \cite{rs} for details about the proof. In particular, we get that

$$\lim_{R\to +\infty}\lim_{k\to +\infty}\int_{B_{R\mk}(\xk)}\Vk e^{4\uk}\, dx=16\pi^2.$$

\medskip\noindent{\bf Step \ref{sec:est}.3:} We claim that there exists $N\in\nn^\star$, there exist $(x_{1,k})_{k\in\nn},...,(x_{N,k})_{k\in\nn}$ such that 
\bequa\label{est:pt:1}
(\inf_{i\in\{1,...,N\}} |x-\xik|) e^{\uk(x)}\leq C(\omega)
\eequa
for all $x\in\omega$ and all $k\in\nn$. Here $x_{1,k}:=\xk$.

\medskip\noindent{\it Proof of the claim:} If there exists $C(\omega)>0$ such that $ |x-\xk| e^{\uk(x)}\leq C(\omega)$ for all $k\in\nn$ and all $x\in\omega$, then we are done. Otherwise, let $\yk\in\overline{\omega}$ such that
\bequa\label{lim:construc:pts}
\sup_{x\in\omega} |x-\xk| e^{\uk(x)}= |\yk-\xk| e^{\uk(\yk)}\to +\infty
\eequa
when $k\to +\infty$. We define
$$\hat{u}_k(x):=\uk(\yk+\nu_k x)-\uk(\yk)$$
for all $x\in \nu_k^{-1}(\omega-\yk)$, where $\nu_k=e^{-\uk(yk)}$ for all $k\in\nn$. It follows from \eqref{lim:construc:pts} that $\hat{u}_k$ is bounded from above uniformly locally on $\rn$ independantly of $k$. We proceed as in Steps \ref{sec:est}.1 and \ref{sec:est}.2 and prove that $\hat{u}_k$ converges to $U_0$ in $C^3_{loc}(\rn)$, and that these two rescaled functions do not interact one with the other. We then add another level of energy $16\pi^2$. If \eqref{est:pt:1} holds with $x_{1,k}=\xk$ and $x_{2,k}=\yk$, then we are done. Otherwise, the process goes on and must cease, because when we have constructed $N$ points, we have that the energy $16\pi^2 N$ and with \eqref{hyp:lambda} we must have $16\pi^2N\leq \Lambda$. We refer to \cite{rs} for the details.\hfill$\Box$

\medskip\noindent{\bf Step \ref{sec:est}.4:} We claim that 
$$\lim_{k\to +\infty}\uk= -\infty$$
uniformly on every compact subset of $\Omega\setminus S_0$. 

\medskip\noindent{\it Proof of the claim:} We prove the claim by contradiction and assume that the conclusion is false. It then follows that point (A1) of Proposition \ref{bound:out0} holds, and then that $\uk$ is uniformly bounded in $C^3_{loc}(\Omega\setminus S_0)$. We let $x_0\in S_0$ and $\delta>0$ such that $B_{2\delta}(x_0)\subset\Omega$ and $B_{2\delta}(x_0)\cap S_0=\{x_0\}$. We let $\xk\in\Omega$ such that $\uk(\xk)=\sup_{B_{\delta}(x_0)}\uk$ and we define $\vk$ and $\mk$ as in \eqref{def:vk:proof}. As in the proof of Proposition \ref{estimate:pointwise1}, we get that $\lim_{k\to +\infty}\xk=x_0\in\omega\cap S_0$ and that \eqref{lim:vk:ellip} holds. We let $H_\delta$ be the Green's function for $\Delta^2$ in $B_{\delta}(x_0)$ with Navier condition on the boundary, that is for any $x\in B_\delta(x_0)$, we have that $\Delta^2 H_\delta(x,\cdot)=\delta_x$ in ${\mathcal D}'(B_\delta(x_0))$ and $H_\delta(x,\cdot)=\Delta H_{\delta}(x,\cdot)=0$ on $\partial B_{\delta}(x_0)$. For any $x\in B_{\delta}(x_0)\setminus\{x_0\}$, we then have that
$$\uk(x)=\int_{B_\delta(x_0)}H_\delta(x,y)\Vk(y) e^{4\uk(y)}\, dy+\varphi_k(x)$$ 
where $\Delta^2\varphi_k=0$, $\varphi_k(y)=\uk(y)$ and $\Delta\varphi_k(y)=\Delta\uk(y)$ for $y\in \partial B_{\delta}(x_0)$. It follows from point (A1) of Proposition \ref{bound:out0} and the comparison principle that $\varphi_k$ is uniformly bounded when $k\to +\infty$. Since $H_\delta>0$, we get with (\ref{lim:Vk}), $(E)$, a change of variable and \eqref{eq:vk:new} that
\beq
\uk(x) &\geq &  \int_{B_{R\mk}(\xk)}H_\delta(x,y) \Vk(y) e^{4\uk(y)}\, dy-C\\
&\geq & \int_{B_{R}(0)}H_\delta(x,\xk+\mk x) \tVk(y)e^{4\vk(y)}\, dy-C\\
&\geq & \int_{B_{R}(0)}H_\delta(x,x_0)\lim_{k\to +\infty}\left(\tVk(y) e^{4\vk(y)}\right)\, dy
\eeq
With standard properties of $H_\delta$, we get that $H_\delta(x,x_0)\geq \frac{1}{8\pi^2}\ln\frac{1}{|x-x_0|}-C$ for $x\in B_{\delta/2}(x_0)$. We then get with \eqref{lim:vk:ellip} that
$$\uk(x)\geq 2\ln\frac{1}{|x-x_0|}-C'$$
for $x\in B_{\delta/2}(x_0)$, $x\neq x_0$ and $k$ large depending on a lower bound on $|x-x_0|$. We then get that for any $0<\alpha<\beta$ small,
$$\Lambda\geq \int_{B_\beta(x_0)\setminus B_\alpha(x_0)} \Vk e^{4\uk}\, dx\geq C\int_{B_\beta(x_0)\setminus B_\alpha(x_0)}\frac{1}{|x-x_0|^8}\, dx,$$
for $k$ large depending on $\alpha$. We then get a contradiction by letting $\alpha\to 0$. Then Case (A1) of Proposition \ref{bound:out0} does not hold and Case (A2) holds. We then get that $\lim_{k\to +\infty}\uk=-\infty$ on compact subsets of $\Omega\setminus S_0$ when $k\to +\infty$. This proves the claim.\hfill$\Box$

\medskip\noindent{\bf Step \ref{sec:est}.5:} We claim that for any $\omega\subset\subset\Omega$, there exists $C(\omega)>0$ such that
\bequa\label{ineq:delta}
(\inf_{i\in\{1,...,N\}}|x-\xik|)^2|\Delta\uk(x)|\leq C(\omega)
\eequa
for all $x\in\omega$ and all $k\in\nn$.

\medskip\noindent{\it Proof of the claim:} We let $x_0\in S_0$ and let $\delta>0$ such that $B_{3\delta}(x_0)\subset\Omega$ and $B_{3\delta}(x_0)\cap S_0=\{x_0\}$. We denote $H_\delta$ the Green's function for $\Delta$ on $B_{2\delta}(x_0)$ with Dirichlet boundary condition. It follows from Green's representation formula that
\bequa\label{id:green}
\Delta\uk(x)=\int_{B_{2\delta}(x_0)}H_\delta(x,y)\Delta^2\uk(y)\, dy + \psi_k(x)\eequa
for all $x\in B_{2\delta}(x_0)$. In this expression, $\psi_k$ is such that $\Delta\psi_k=0$ in $B_{2\delta}(x_0)$ and $\psi_k(x)=\Delta\uk(x)$ on $\partial B_{2\delta}(x_0)$. It follows from Proposition \ref{bound:out0} and the comparison principle that there exists $C_\delta>0$ such that
\bequa\label{ineq:psi:1}
|\psi_k(x)|\leq C_\delta
\eequa
for all $x\in B_{2\delta}(x_0)$. We consider a sequence $(\yk)_{k\in\nn}\in B_\delta(x_0)$ that converges. We assume that $\lim_{k\to +\infty}\yk=x_0$ when $k\to +\infty$. With standard properties of the Green's function, \eqref{id:green} and \eqref{ineq:psi:1}, we get that there exists $C>0$ such that
$$|\Delta\uk(\yk)|\leq C\int_{B_{2\delta}(x_0)}\frac{e^{4\uk(y)}}{|\yk-y|^2}\, dy+C.$$
We let $\Rk(x)=\inf_{i\in\{1,...,N\}}|x-\xik|$ for all $x\in\Omega$, we let $\theta_{i,k}=\frac{\yk-\xik}{|\yk-\xik|}$ and $\Omega_{i,k}=\{y\in B_{2\delta}(x_0)/ \Rk(y)=|y-\xik|\}$.  With \eqref{hyp:lambda} and the pointwise estimate \eqref{est:pt:1}, we then get that
\beq
|\Delta\uk(\yk)| &\leq&  C \int_{B_{2\delta}(x_0)\setminus \cup B_{\frac{|\yk-\xik|}{2}}(\xik)}+C\int_{\cup B_{\frac{|\yk-\xik|}{2}}(\xik)}+C\\
&\leq &  C\sum_{i=1}^N \int_{\Omega_{i,k}\setminus B_{\frac{|\yk-\xik|}{2}}(\xik)}+C\sum_{i=1}^N \int_{B_{\frac{|\yk-\xik|}{2}}(\xik)} +C\\
&\leq & C\sum_{i=1}^N \int_{B_{2\delta}(x_0)\setminus B_{\frac{|\yk-\xik|}{2}}(\xik)}\frac{1}{|\yk-y|^2 |\xik-y|^4}\, dy\\
& & +C\sum_{i=1}^N \int_{B_{\frac{|\yk-\xik|}{2}}(\xik)}\frac{e^{4\uk(y)}}{\Rk(\yk)^2}\, dy+ C\\
&\leq &  C\sum_{i=1}^N \int_{\rn\setminus B_{\frac{1}{2}}(0)}\frac{1}{\Rk(\yk)^2 |\theta_{i,k}-z|^2 |z|^4}\, dz +\frac{C}{\Rk(\yk)}+ C
\eeq
for all $k\in\nn$ large enough, and then
\bequa\label{bnd:sup:delta}
\Rk(\yk)^2|\Delta\uk(\yk)|=O(1)
\eequa
when $k\to +\infty$ in case $\lim_{k\to +\infty}\yk=x_0$. When $\lim_{k\to +\infty}\yk\neq x_0$, inequality \eqref{bnd:sup:delta} is a consequence of Proposition \ref{bound:out0}. Since the sequence $\yk$ is arbitrary, this proves \eqref{ineq:delta} on $B_{\delta}(x_0)$. As easily checked, \eqref{ineq:delta} follows from this estimate taken in the neighborhood of each of the points in $S_0$ and Proposition \ref{bound:out0}.\hfill$\Box$

\medskip\noindent Proposition \eqref{estimate:pointwise1} is a consequence of Steps \ref{sec:est}.1 to \ref{sec:est}.5.

\section{Blow-Up analysis}\label{sec:blowup}
The proof of Theorem \ref{th:intro} goes through an induction that will use the following proposition. The paper of Li-Shafrir \cite{lishafrir} was a source of inspiration. 

\begin{prop}\label{prop:funda} 
Let $x_0\in\rn$, $\delta>0$ and $\Lambda>0$. We let $\Vk\in C^0(B_{4\delta}(x_0))$ such that $\lim_{k\to +\infty}\Vk=1$ in $C^0(B_{4\delta}(x_0))$. We let $\uk\in C^4(B_{4\delta}(x_0))$ such that
\bequa\label{eq:ue:p}
\Delta^2\uk=\Vk e^{4\uk}
\eequa
in $B_{4\delta}(x_0)$. We assume that
\bequa\label{bnd:nrj:p}
\int_{B_{4\delta}(x_0)}e^{4\uk}\, dx\leq \Lambda
\eequa
for all $k\in\nn$. We let $\rhok\geq 0$ such that $\lim_{k\to +\infty}\rhok=0$. We assume that there exists $(\xk=x_{1,k})_{k\in\nn},...,(x_{N,k})_{k\in\nn}\in B_{4\delta}(x_0)$ such that for any $i\in \{1,...,N\}$, we have that
\bequa\label{lim:xie:p}
\lim_{k\to +\infty}\xik=x_0\hbox{ and }\lim_{k\to +\infty}\uk(\xik)=+\infty.
\eequa
Moreover, we assume that there exists $C>0$ such that
\bequa\label{ineq:wpe:p}
\inf_{i\in \{1,...,N\}}|x-\xik|e^{\uk(x)}\leq C\hbox{ and }\inf_{i\in \{1,...,N\}}|x-\xik|^2|\Delta\uk(x)|\leq C
\eequa
for all $k\in\nn$ and all $x\in B_{2\delta}(\xk)\setminus \ovb_{\rhok}(\xk)$. We assume that
\bequa\label{lim:xie:infty:p}
\lim_{k\to +\infty}\frac{|\xik-\xjk|}{\mik}=+\infty
\eequa
for all $i\neq j$, $i,j\in \{1,...,N\}$. In this expression, we have let $\mik=e^{-\uk(\xik)}$. We assume that
\bequa\label{cv:resc:ue:p}
\lim_{k\to +\infty}(\uk(\xik+\mik x)-\uk(\xik))=\ln\frac{\sqrt{96}}{\sqrt{96}+|x|^2}
\eequa
for all $x\in\rn$, and that this convergence holds in $C^3_{loc}(\rn)$. We let $(\rk)_{k\in\nn}$ such that $\rk>0$ for all $k\in\nn$ and that $\lim_{k\to +\infty}\rk=r\in [0,\delta]$. We let
\bequa\label{def:I:p}
I:=\left\{i\in\{2,...,N\}\,/\,\frac{\xik-\xk}{\rk}=O(1)\hbox{ when }k\to +\infty\right\}.\eequa
Note that $I$ may be empty. We let $\txi=\lim_{k\to +\infty}\frac{\xik-\xk}{\rk}$ for $i\in I$. We assume that $\txi\neq 0$ for all $i\in I$ and that
\bequa\label{lim:rhoe:p}
\rhok=o(\rk)
\eequa
and that $\mk=\mu_{1,k}=o(\rk)$ when $k\to +\infty$. We let $\nu,R$ such that
\bequa\label{choice:nu:p}
0<\nu< \frac{1}{10}\min\left\{\{|\txi|/\, i\in I\}\cup\{|\txi-\txj|/\, i,j\in I,\, \txi\neq\txj\}\right\}
\eequa
and
\bequa\label{choice:R:p}
3\max\left\{|\txi|/\, i\in I\right\}<R<\frac{\delta}{2r}.
\eequa
In case $r=0$, we let $\frac{\delta}{2r}=+\infty$. We let
$$D_k:= B_{R\rk}(\xk)\setminus\bigcup_{i\in I}\ovb_{\nu\rk}(\xik).$$
Then, if $\mk=o(\rhok)$, we have that
$$\lim_{k\to +\infty}\int_{D_k\setminus \ovb_{2\rhok}(\xk)}e^{4\uk(x)}\, dx=0.$$If $\rhok=O(\mk)$, we have that
$$\lim_{\tilde{R}\to +\infty}\lim_{k\to +\infty}\int_{D_k\setminus \ovb_{\tilde{R}\mk}(\xk)}e^{4\uk(x)}\, dx=0.$$
\end{prop}
This section is devoted to the proof of the proposition. Up to relabelling the $\txi$'s, we assume that there exists $\phi:\{1,...,l\}\to \{1,...,N\}$ such that $\txpi\neq \txpj$ for all $i\neq j$, and
\bequa\label{relabel}
\left\{\txi/\, i\in I\right\}=\left\{\txpi/\, i\in \{1,..,l\}\right\}, \hbox{ and }I=\{2,...,\phi(l)\}.
\eequa
Moreover, we assume that $\phi$ is increasing and $\txj=\txpi$ for all $j$ such that $\phi(i)\leq j<\phi(i+1)$. Note that $1\not\in I$ and that $\phi(1)\neq 1$. For all $i\in \{1,...,N\}$, we let
\bequa\label{def:txie:p}
\txik:=\frac{\xik-\xk}{\rk}.
\eequa

\medskip\noindent{\bf Step \ref{sec:blowup}.1 (Rescaling):} We let 
$$\Omega_k=\left(B_{3R}(0)\setminus\bigcup_{i=1}^l\ovb_{\nu}(\txpik)\right)\setminus\ovb_{\frac{\rhok}{\rk}}(0).$$
With the choice (\ref{choice:R:p}) of $R$, we have that $\xk+\rk x\in B_{2\delta}(\xk)\subset B_{4\delta}(x_0)$ for all $x\in\Omega_k$. We let $x\in\Omega_k$, and $j\in \{1,...,N\}$. We distinguish three cases:

\smallskip\noindent{\it Case \ref{sec:blowup}.1.1:} We assume that $j\in I$. We let $i\in \{1,...,l\}$ such that $\phi(i)\leq j<\phi(i+1)$. Then, with (\ref{choice:nu:p}), (\ref{choice:R:p}), the definition \eqref{def:txie:p} and the choice of the numbering of the $\txj$'s, we have that
\beqn
|\xk+\rk x-\xjk|&=&\rk |x-\txjk|\geq \rk\left(|x-\txpik|-|\txpik-\txjk|\right)\nonumber\\
&\geq& \rk (\nu+o(1))\geq \rk\frac{\nu}{2}\geq \rk\frac{\nu}{6R} |x|.\label{case1:sec3}
\eeqn

\smallskip\noindent{\it Case \ref{sec:blowup}.1.2:} We assume that $j\in\{2,...,N\}$ is such that $j\not\in I$. Then with (\ref{choice:R:p}), the definition \eqref{def:txie:p} and the definition (\ref{def:I:p}) of $I$, we get that
\bequa\label{case2:sec3}
|\xk+\rk x-\xjk|=\rk |x-\txjk|\geq \rk\left(|\txjk|-3R\right)\geq \rk\geq\frac{\rk}{3R} |x|.
\eequa

\smallskip\noindent{\it Case \ref{sec:blowup}.1.3:} If $j=1$, we get that
\bequa\label{case3:sec3}
|\xk+\rk x-\xk|=\rk |x|>\rhok.
\eequa

\smallskip\noindent It follows from \eqref{case1:sec3}-\eqref{case3:sec3} that 
\bequa\label{ineq:proof:case1}
\inf_{i\in \{1,...,N\}}|\xk+\rk x-\xik|\geq C(\nu,R)\rk |x|
\eequa
$$\hbox{ and }\xk+\rk x\in B_{2\delta}(\xk)\setminus \ovb_{\rhok}(\xk).$$
for all $x\in \Omega_k$. For $x\in B_{3R}(0)$, we let
\bequa\label{def:tue:p}
\tuk(x):=\uk(\xk+\rk x)+\ln\rk.
\eequa
It follows from \eqref{ineq:proof:case1}, (\ref{eq:ue:p}) and (\ref{ineq:wpe:p}) that there exists $C>0$ such that 
\bequa\label{eq:tue:p}
\Delta^2\tuk= \tVk e^{4\tuk}\hbox{ in }B_{3R}(0)
\eequa
and
\bequa\label{ineq:wpe:tue:p}
|x|e^{\tuk(x)}\leq C\hbox{ and }|x|^2|\Delta\tuk(x)|\leq C
\eequa
for all $x\in\Omega_k$. Here, we let $\tVk(x):=\Vk(\xk+\rk x)$ for all $x\in B_{3R}(0)$ and all $k\in\nn$.

\medskip\noindent{\bf Step \ref{sec:blowup}.2 (Harnack inequality):} We claim that there exists $C=C(\nu,R)$, there exists $\beta=\beta(\nu,R)>0$ such that
\bequa\label{ineq:harnack}
\beta\sup_{\partial\left(B_r(0)\setminus\bigcup_{i=1}^l\ovb_{2\nu}(\txpik)\right)}\tuk\leq \inf_{\partial\left(B_r(0)\setminus\bigcup_{i=1}^l\ovb_{2\nu}(\txpik)\right)}\tuk+(1-\beta)\ln r+C
\eequa
for all $r>0$ such that
$$\frac{3\rhok}{\rk}\leq r\leq 2R.$$

\medskip\noindent{\it Proof of the claim:} We let $\sk>0$ such that $\frac{3\rhok}{\rk}\leq \sk\leq 2R$. Up to a subsequence, we assume that $\lim_{k\to +\infty}\sk=s\geq 0$. We distinguish two cases:

\smallskip\noindent{\it Case \ref{sec:blowup}.2.1:} We assume that
\bequa\label{hyp:s:1}
0\leq s<\frac{4}{5}\nu.
\eequa
With (\ref{choice:nu:p}) and (\ref{choice:R:p}), we get that
$$B_{\frac{5}{4}}(0)\setminus\ovb_{\frac{1}{2}}(0)\subset\frac{\Omega_k}{\sk}.$$
For any $x\in B_{\frac{5}{4}}(0)\setminus\ovb_{\frac{1}{2}}(0)$, we define
\bequa\label{def:Ue:p}
\Uk(x)=\tuk(\sk x)+\ln\sk.
\eequa
It follows from (\ref{ineq:wpe:tue:p}) that there exists $C>0$ such that 
\bequa\label{estim:Ue:p}
\Uk(x)\leq C\hbox{ and }|\Delta\Uk(x)|\leq C
\eequa
for all $k\in\nn$ and all $x\in B_{\frac{5}{4}}(0)\setminus\ovb_{\frac{1}{2}}(0)$. It then follows from the Harnack inequality that there exists $\beta, C>0$ such that
\bequa\label{ineq:har:1}
\beta\sup_{\partial B_1(0)}\Uk\leq \inf_{\partial B_1(0)}\Uk+C
\eequa
for all $k\in\nn$. Coming back to $\tuk$ with (\ref{def:Ue:p}), using the assumption \eqref{hyp:s:1}, (\ref{choice:nu:p}) and (\ref{choice:R:p}) we get that
$$\partial\left(B_{\sk}(0)\setminus\bigcup_{i=1}^l\ovb_{2\nu}(\txpik)\right)=\partial B_{\sk}(0)$$
and then \eqref{ineq:harnack} follows from \eqref{ineq:har:1}. This ends Case \ref{sec:blowup}.2.1.

\medskip\noindent{\it Case \ref{sec:blowup}.2.2:} We assume that
\bequa\label{ineq:case:2}
\frac{4}{5}\nu\leq s\leq 2R.
\eequa
We let 
$${\mathcal A}=\left(B_{3R}(0)\setminus\bigcup_{i=1}^l\ovb_{\frac{5}{4}\nu}(\txpi)\right)\setminus\ovb_{\frac{\nu}{5}}(0).$$
It follows from (\ref{choice:nu:p}) and (\ref{lim:rhoe:p}) that 
$${\mathcal A}\subset \Omega_k$$
for $k>0$ large enough. Moreover, it follows from (\ref{choice:nu:p}) and (\ref{choice:R:p}) that the balls $\ovb_{\frac{\nu}{5}}(0)$, $\ovb_{\frac{5}{4}\nu}(\txpi)$, $(i\in \{1,...,l\})$ are disjoint and contained in $B_{2R}(0)$. We then get that ${\mathcal A}$ is connected. 

\noindent It follows from (\ref{ineq:wpe:tue:p}) that there exists $C>0$ such that
$$\tuk(x)\leq C\hbox{ and }|\Delta\tuk(x)|\leq C$$
for all $x$ in a neighborhood of ${\mathcal A}$. With Harnack's inequality, we get that there exists $\beta,C>0$ such that


$$\beta\sup_{{\mathcal A}}\tuk\leq \inf_{{\mathcal A}}\tuk+C$$
for all $k\in\nn$. With (\ref{choice:nu:p}), (\ref{choice:R:p}) and $\frac{4}{5}\nu\leq s\leq 2R$, we get that 
$$\partial\left(B_{\sk}(0)\setminus\bigcup_{i=1}^l\ovb_{2\nu}(\txpik)\right)\subset{\mathcal A}.$$
With \eqref{ineq:case:2}, we get that there exists $\beta=\beta(\nu,R)>0$, $C=C(\nu,R)>0$ such that
$$\beta\sup_{\partial\left(B_{\sk}(0)\setminus\bigcup_{i=1}^l\ovb_{2\nu}(\txpik)\right)}\tuk\leq \inf_{\partial\left(B_{\sk}(0)\setminus\bigcup_{i=1}^l\ovb_{2\nu}(\txpik)\right)}\tuk+(1-\beta)\ln \sk+C$$
for all $k\in\nn$. This ends Case \ref{sec:blowup}.2.2, and the proof of the claim is complete.\hfill$\Box$

\medskip\noindent{\bf Step  \ref{sec:blowup}.3 (Upper bound):} We claim that there exists $\theta>-1$, there exists $R_0>0$ such that
\bequa\label{ineq:upper}
\sup_{\partial\left(B_{\sk}(0)\setminus\bigcup_{i=1}^l\ovb_{2\nu}(\txpik)\right)}\tuk\leq -\left(1+\frac{1+\theta}{\beta}\right)\ln\sk-\frac{1+\theta}{\beta}\ln\frac{\rk}{\mk}+C
\eequa
for all $k\in\nn$ where $\sk>0$ is such that
\bequa\label{choice:sk:1}
\sk\in \left[\frac{3\rhok}{\rk},2R\right]\hbox{ if
}\mk=o(\rhok)
\eequa
and
\bequa\label{choice:sk:2}
\sk\in\left[\frac{R_0\mk}{\rk},2R\right]\hbox{ if }\rhok=O(\mk).
\eequa

\medskip\noindent{\it Proof of the claim:} We let $\Uk$ defined as in (\ref{def:Ue:p}) on $B_{\frac{3R}{\sk}}(0)$. We assume that 
\bequa\label{hyp:s:nu}
0\leq s<8\nu.
\eequa
Let $H_k$ be the Green's function of $\Delta^2$ on 
$${\mathcal D}_{k}:=B_{1}(0)\setminus\bigcup_{i=1}^l\ovb_{\frac{2\nu}{\sk}}\left(\frac{\txpik}{\sk}\right)=B_1(0)$$
with Navier condition on the boundary, that is for any $x\in{\mathcal D}_k$, we have that $\Delta^2H_{\delta}(x,\cdot)=\delta_x$ in the distribution sense, $H_\delta(x,\cdot)=\Delta H_\delta(x,\cdot)=0$ on $\partial\Omega$. Note that the preceding inequality is a consequence of (\ref{choice:nu:p}), (\ref{choice:R:p}) and \eqref{hyp:s:nu}. With (\ref{choice:nu:p}) and (\ref{choice:R:p}), we get that $0\in {\mathcal D}_{k}$. It follows from Green's representation formula that
\bequa\label{eq:Uk}
\Uk(0)=\int_{{\mathcal D}_{k}}H_k(0,y)\Delta^2\Uk(y)\, dy+\varphi_k(0)+\psi_k(0)\eequa
where 
\bequa\label{def:phi:psi}
\left\{\begin{array}{ll}
\Delta\varphi_k=0 & \hbox{ in } {\mathcal D}_{k}\\
\varphi_k=\Uk & \hbox{ on }\partial {\mathcal D}_{k}
\end{array}\right\}\hbox{ and }
\left\{\begin{array}{ll}
\Delta^2\psi_k=0 & \hbox{ in } {\mathcal D}_{k}\\
\Delta\psi_k=\Delta\Uk & \hbox{ on }\partial {\mathcal D}_{k}\\
\psi_k=0 & \hbox{ on }\partial {\mathcal D}_{k}
\end{array}\right\}.
\eequa
It follows from \eqref{lim:rhoe:p}, (\ref{choice:nu:p}) and (\ref{choice:R:p}) that
$$\partial {\mathcal D}_{k}=\partial B_1(0)\subset \frac{\Omega_k}{\sk}.$$
We then get with (\ref{estim:Ue:p}), the maximum principle and
(\ref{def:phi:psi}) that there exists $C>0$ such that
\bequa\label{ineq:psi:2}
|\psi_k(0)|\leq C
\eequa
for all $k\in\nn$. It follows from the comparison principle and \eqref{def:phi:psi} that
\bequa\label{ineq:phi:2}
\varphi_k(0)\geq \inf_{\partial {\mathcal D}_{k}}\Uk.
\eequa
We let $\tilde{R}>0$. Moreover, with \eqref{choice:sk:1}, \eqref{choice:sk:1} we get that
\bequa\label{incl:b}
B_{\frac{\tilde{R}\mk}{\sk\rk}}(0)\subset \partial B_{1/2}(0)\subset{\mathcal D}_k
\eequa
with $R_0>2\tilde{R}$. Here, we have let $\mu_k=\mu_{1,k}=e^{-\uk(\xk)}$. Noting that $H_k\geq 0$, we get with \eqref{eq:Uk}, \eqref{ineq:psi:2}, \eqref{ineq:phi:2} and \eqref{incl:b} that
$$\Uk(0)\geq \int_{B_{\frac{\tilde{R}\mk}{\sk\rk}}(0)}H_k(0,y)\Delta^2 \Uk(y)\, dy+\inf_{\partial {\mathcal D}_{k}}\Uk-C.$$
It follows from standard elliptic estimates that there exists $C>0$ such that
$$H_k(0,y)\geq \frac{1}{8\pi^2}\ln\frac{1}{|y|}-C$$
for $y\in B_{1/2}(0)\setminus\{0\}$. We then get that
$$\Uk(0)\geq \int_{B_{\frac{\tilde{R}\mk}{\sk\rk}}(0)}\left(\frac{1}{8\pi^2}\ln\frac{1}{|y|}-C\right) \Delta^2\Uk(y)\, dy+\inf_{\partial {\mathcal D}_{k}}\Uk-C.$$
With the change of variable $y=\frac{\mk}{\sk\rk}z$ (where $\mu_k=\mu_{1,k}=e^{-\uk(\xk)}$) and coming back to the definitions (\ref{def:tue:p}) and (\ref{def:Ue:p}), we get that
\beq
\ln\frac{\rk\sk}{\mk}&\geq &\int_{B_{\tilde{R}}(0)}\left(\frac{1}{8\pi^2}\ln\frac{\sk\rk}{\mk}+\frac{1}{8\pi^2}\ln\frac{1}{|z|}-C\right) \tVk(x) e^{4(\uk(\xk+\mk z)-\uk(\xk))}\, dz\\
&&+\inf_{\partial {\mathcal D}_{k}}\Uk-C.
\eeq
With (\ref{cv:resc:ue:p}), we then get that

$$C(\tilde{R})\geq\left(1+\frac{\theta_k(\tilde{R})}{8\pi^2}\right)
\ln\frac{\sk\rk}{\mk}+\inf_{\partial {\mathcal D}_{k}}\Uk,$$
with $\lim_{\tilde{R}\to +\infty}\lim_{k\to
  +\infty}\theta_k(\tilde{R})=0$. Choosing $\tilde{R}$ large enough, and then choosing $R_0>2\tilde{R}$ large, we get that there exists $\theta>-1$ such that
$$C\geq\left(1+\theta\right) \ln\frac{\sk\rk}{\mk}+\inf_{\partial {\mathcal D}_{k}}\Uk,$$
for all $k\in\nn$. Coming back to $\tuk$ and using \eqref{ineq:harnack}, we get the inequality of the Lemma. This ends the proof of the claim when \eqref{hyp:s:nu} holds. In case \eqref{hyp:s:nu} does not hold, the claim follows from the case $\sk=7\nu$ and the Harnack inequality \eqref{ineq:harnack}.\hfill$\Box$

\medskip\noindent{\bf Step \ref{sec:blowup}.4 (Proof of Proposition \ref{prop:funda}):} We let $\yk\in B_{2R}(0)\setminus\bigcup_{i=1}^l\ovb_{3\nu}(\txpik)$ such that
$$|\yk|\geq \frac{3\rhok}{\rk}\hbox{ if }\mk=o(\rhok)\hbox{ or }|\yk|\geq \frac{R_0\mk}{\rk}\hbox{ if }\rhok=O(\mk),$$
where $R_0$ is as in \eqref{ineq:upper}. We let $\sk=|\yk|$, so that
$$\yk\in \partial \left(B_{\sk}(0)\setminus\bigcup_{i=1}^l\ovb_{2\nu}(\txpik)\right).$$
It follows from \eqref{ineq:upper} that
\bequa\label{ineq:proof:u}
\tuk(\yk)\leq
-\left(1+\frac{1+\theta}{\beta}\right)\ln|\yk|-\frac{1+\theta}{\beta}\ln\frac{\rk}{\mk}+C.
\eequa
We distinguish two cases:

\smallskip\noindent{\it Case \ref{sec:blowup}.4.1:} We assume that
$\mk=o(\rhok)$. We then get with \eqref{ineq:proof:u} that
\beq
&&\int_{\left( B_{2R}(0)\setminus\bigcup_{i=1}^l\ovb_{2\nu}(\txpik)\right)\setminus \ovb_{\frac{3\rhok}{\rk}}(0)}e^{4\tuk(y)}\, dy\\
&&\leq C\int_{B_{2R}(0)\setminus \ovb_{\frac{3\rhok}{\rk}}(0)}\left(\frac{\mk}{\rk}\right)^{4\frac{1+\theta}{\beta}}\frac{1}{|y|^{4+4\frac{1+\theta}{\beta}}}\, dy\leq C\left(\frac{\mk}{\rhok}\right)^{4\frac{1+\theta}{\beta}}=o(1)
\eeq
when $k\to +\infty$. Coming back to the definition of $\tuk$ and the relabelling \eqref{relabel}, this proves Proposition \ref{prop:funda} in Case \ref{sec:blowup}.4.1.

\smallskip\noindent{\it Case \ref{sec:blowup}.4.2:} We assume that
$\rhok=O(\mk)$ when $k\to +\infty$. We take $\tilde{R}>R_0$. We then get
with \eqref{ineq:proof:u} that
\beq
&&\int_{\left( B_{2R}(0)\setminus\bigcup_{i=1}^l\ovb_{2\nu}(\txpik)\right)\setminus \ovb_{\frac{\tilde{R}\mk}{\rk}}(0)}e^{4\tuk(y)}\, dy\\
&&\leq C\int_{B_{2R}(0)\setminus \ovb_{\frac{\tilde{R}\mk}{\rk}}(0)}\left(\frac{\mk}{\rk}\right)^{4\frac{1+\theta}{\beta}}\frac{1}{|y|^{4+4\frac{1+\theta}{\beta}}}\, dy\leq \frac{C}{\tilde{R}^{4\frac{1+\theta}{\beta}}}
\eeq
for all $k\in\nn$. Coming back to the definition of $\tuk$ and the relabelling \eqref{relabel}, this proves Proposition \ref{prop:funda} in Case \ref{sec:blowup}.4.2.

\section{Proof of Theorem \ref{th:intro}}\label{sec:conc}
We prove Theorem \ref{th:intro} by induction. We let $N\in \mathbb{N}^\star$. We say that ${\bf (H}_{\bf N}{\bf )}$ holds if the following Proposition holds:

\noindent{\bf Proposition ${\bf (H}_{\bf N}{\bf )}$: }{\it  Let $x_0\in\rn$, $\delta>0$ and $\lambda>0$. Let $\uk\in C^4(B_{4\delta}(x_0))$ such that
\bequa\label{eq:rec}
\Delta^2\uk=\Vk e^{4\uk}
\eequa
in $B_{4\delta}(x_0)$ and
$$\int_{B_{4\delta}(x_0)}e^{4\uk}\, dx\leq \Lambda.$$
We assume that there exists $1\leq K\leq N$, $\xk=x_{1,k},...,x_{K,k}\in B_{4\delta}(x_0)$ such that for any $i\in \{1,...,K\}$, we have that
$$\lim_{k\to +\infty}\xik=x_0\hbox{ and }\lim_{k\to +\infty}\uk(\xik)=+\infty.$$
Moreover, we assume that there exists $C>0$ such that
\bequa\label{estim:1}
\inf_{i\in \{1,...,K\}}|x-\xik|e^{\uk(x)}\leq C\hbox{ and }\inf_{i\in \{1,...,K\}}|x-\xik|^2|\Delta\uk(x)|\leq C
\eequa
for all $k\in\nn$ and all $x\in B_{2\delta}(\xk)$. We assume that
\bequa\label{estim:2}
\lim_{k\to +\infty}\frac{|\xik-\xjk|}{\mik}=+\infty
\eequa
for all $i\neq j$, $i,j\in \{1,...,K\}$. In this expression, we have let $\mik=e^{-\uk(\xik)}$. We assume that
\bequa\label{estim:3}
\lim_{k\to +\infty}(\uk(\xik+\mik x)-\uk(\xik))=\ln\frac{\sqrt{96}}{\sqrt{96}+|x|^2}
\eequa
for all $x\in\rn$, and that this convergence holds in $C^3_{loc}(\rn)$. Then, we have that
$$\int_{B_\delta(x_0)} \Vk e^{4\uk(x)}\, dx=16\pi^2K+o(1)$$
when $k\to +\infty.$}

We prove by induction that ${\bf (H}_{\bf N}{\bf )}$ holds for all $N\geq 1$.

\medskip\noindent{\bf Step \ref{sec:conc}.1 (Proof of ${\bf (H}_{\bf 1}{\bf )}$):} We claim that ${\bf (H}_{\bf 1}{\bf )}$ holds. We prove the claim. We apply Proposition \ref{prop:funda} with $\rk=\delta$ and $\rhok=0$. We then get that
\bequa\label{lim:rec:init}
\lim_{R\to +\infty}\lim_{k\to +\infty}\int_{B_{\frac{\delta}{2}}(x_0)\setminus B_{R\mk}(\xk)}\Vk e^{4\uk(x)}\, dx=0.
\eequa
Plugging (\ref{estim:3}), Proposition \ref{estimate:pointwise1} and \eqref{lim:rec:init} together yields
$$\int_{B_\delta(x_0)}\Vk e^{4\uk(x)}\, dx=16\pi^2+o(1)$$
when $k\to +\infty$. This proves the claim, and therefore ${\bf (H}_{\bf 1}{\bf )}$.

\medskip\noindent{\bf Step \ref{sec:conc}.2 (Induction):} We let $N\geq 2$. We assume that ${\bf (H}_{\bf N-1}{\bf )}$ holds. We let $(\uk)_{k\in\nn}\in C^4(B_{4\delta}(x_0))$. We assume that $u_k$ verifies the hypothesis of ${\bf (H}_{\bf N}{\bf )}$. Clearly we can assume that $K=N$ in the statement of ${\bf (H}_{\bf N}{\bf )}$. Up to renumbering, we let
$$r_{1,k}=\inf_{i\neq j}\{|\xik-\xjk|\}=\inf_{i\neq 1}\{|x_{1,k}-\xik|\}.$$
With (\ref{estim:2}), we get that 
$$\lim_{k\to +\infty}\frac{r_{1,k}}{\mu_{1,k}}=+\infty.$$
We let
$$I_1=\left\{i\in\{2,...,N\}/\, \frac{x_{1,k}-\xik}{r_{1,k}}=O(1)\hbox{ as }k\to +\infty\right\}.$$
Note here that $I_1\neq\emptyset$. We define by induction:
\bequa\label{def:rpq}
r_{q+1,k}=\inf\{|x_{1,k}-\xjk|/\, j\not\in \{1\}\cup I_1\cup...\cup
I_q\}
\eequa
\bequa\label{def:Iq}
I_{q+1}=\left\{i\in\{1,...,N\}/\,
  \frac{x_{1,k}-\xik}{r_{q+1,k}}=O(1)\hbox{ as }k\to +\infty,\;
  j\not\in \{1\}\cup I_1\cup...\cup I_q\right\},
\eequa
when these quantities are defined. Since we have a finite number of points, this process must end. We let $q_0\in\nn$ such that $r_{q,k}$ is defined for $q\in \{1,...,q_0\}$ and not afterwards. Moreover, for any $q<q_0$, we have that
$$\lim_{k\to +\infty}\frac{r_{q+1,k}}{r_{q,k}}=+\infty.$$

\medskip\noindent{\it Step \ref{sec:conc}.2.1:} We claim that 
$$\lim_{k\to +\infty}\int_{B_{R r_{1,k}}(\xk)}\Vk e^{4\uk(x)}\, dx=16\pi^2\, \hbox{Card} (\{1\}\cup I_1),$$
where $\xk=x_{1,k}$. We prove the claim. We apply Proposition \ref{prop:funda} with $\uk$, $\rhok=0$ and $\rk=r_{1,k}$. For $R$, $\nu$ and $\phi$ as in the proof of Proposition \ref{prop:funda}, similarly to what was done for the proof of ${\bf (H}_{\bf 1}{\bf )}$  we get that
\bequa\label{lim:holes}
\lim_{k\to +\infty}\int_{B_{R\rk}(\xk)\setminus\bigcup_{i=1}^l\ovb_{\nu\rk}(\xpik)}e^{4\uk(x)}\, dx=16\pi^2.
\eequa
We fix $i\in \{1,...,l\}$. We let
$$\vk(x):=\uk(\xpik+\rk x)+\ln\rk$$
and $\tVk(x):=\Vk(\xpik+\rk x)$ for all $x\in B_R(0)$ and all $k\in\nn$. With (\ref{eq:rec}), we have that
$$\Delta^2\vk= \tVk e^{4\vk}\hbox{ in }B_{4\nu}(0)\hbox{ and }\int_{B_{4\nu}(0)}e^{4\vk}\, dx\leq\Lambda.$$
For any $j$ such that $\phi(i)\leq j<\phi(i+1)$, we let
$$X_{j,k}=\frac{\xjk-\xpik}{\rk}.$$
It follows from the definition of $\phi$ that $\lim_{k\to
  +\infty}X_{j,k}=0$ for all
$j\in\{\phi(i),...,\phi(i+1)-1\}$. Arguing as in Step
\ref{sec:blowup}.1, and letting $U_i:=\{\phi(i),...,\phi(i+1)-1\}$, we get that 





$$\inf_{j\in U_i}|x-X_{j,k}|e^{\vk(x)}\leq C\hbox{ in }\inf_{j\in U_i}|x-X_{j,k}|^2|\Delta\vk(x)|\leq C$$
for all $x\in B_{4\nu}(0)$. For any $j,m\in \{\phi(i),...,\phi(i+1)-1\}$, $j\neq m$, we have with (\ref{estim:2}) that
$$\frac{|X_{j,k}-X_{m,k}|}{e^{-\vk(X_{j,k})}}=\frac{|\xjk-x_{m,k}|}{e^{-\uk(\xjk)}}\to +\infty$$
when $k\to +\infty$. With (\ref{estim:3}), a straightforward computation shows that for any $j\in \{\phi(i),...,\phi(i+1)-1\}$, we have that for any $x\in\rn$,
$$\vk(X_{j,k}+e^{-\vk(X_{j,k})}x)-\vk(X_{j,k})\to \ln\frac{\sqrt{96}}{\sqrt{96}+|x|^2}$$
when $k\to +\infty$. Moreover, this convergence holds in $C^3_{loc}(\rn)$. We then apply the induction hypothesis ${\bf (H}_{\bf N-1}{\bf )}$ with $\vk$ (which has at most $N-1$ concentration points) and we get that 
$$\int_{B_{\nu\rk}(\xpik)} \Vk e^{4\uk(x)}\, dx=\int_{B_{\nu}(0)}\tVk e^{4\vk(x)}\, dx=(\phi(i+1)-\phi(i)) 16\pi^2+o(1)$$
when $k\to +\infty$. Since this inequality is valid for all $i$, we get with (\ref{lim:holes}) that
\beq
\lim_{k\to +\infty}\int_{B_{R\rk}(\xk)} \Vk e^{4\uk(x)}\, dx&=&16\pi^2+\sum_{i=1}^l(\phi(i+1)-\phi(i)) 16\pi^2+o(1)\\
&=&16\pi^2\,\hbox{Card}(\{1\}\cup I_1).
\eeq
This proves the claim, and then Step \ref{sec:conc}.2.1.

\medskip\noindent{\it Step \ref{sec:conc}.2.2:} We let $q<q_0$ and assume that for any $R>0$ large enough, we have that
\bequa\label{lim:step:conc}
\lim_{k\to +\infty}\int_{B_{R r_{q,k}}(\xk)} \Vk e^{4\uk(x)}\,
dx=16\pi^2\,\hbox{Card}(\{1\}\cup I_1\cup ...\cup I_q).
\eequa
We claim that for any $R>0$ large enough, we have that
$$\lim_{k\to +\infty}\int_{B_{R r_{q+1,k}}(\xk)} \Vk e^{4\uk(x)}\, dx=16\pi^2\,\hbox{Card}(\{1\}\cup I_1\cup ...\cup I_{q+1}).$$
We prove the claim. We let 
$$R_0=\max\left\{\frac{|\xk-\xik|}{r_{q,k}}/\, i\in I_q\right\}.$$
We let $\rhok=R_1 r_{q,k}$ with $R_1>2R_0$ and $\rk=r_{q+1,k}$. We let $i\in \{1\}\cup I_1\cup...\cup I_q$ and $x\in B_{2\delta}(\xk)\setminus B_{\rhok}(\xk)$. We assume that $i\in I_p$, $p\leq q$. With the definitions \eqref{def:rpq} and \eqref{def:Iq}, we get that
$$|x-\xik|\geq \frac{|x-\xk|}{2}$$
for all $x\in B_{2\delta}(\xk)\setminus B_{\rhok}(\xk)$ and $i\in \{1\}\cup I_1\cup...\cup I_q$. We then get with (\ref{estim:1}) that there exists $C>0$ such that
$$\inf_{i\not\in I_1\cup...\cup I_q}|x-\xik|e^{\uk(x)}\leq C\hbox{ and
  }\inf_{i\not\in I_1\cup...\cup I_q}|x-\xik|^2|\Delta\uk(x)|\leq C$$
for all $k\in\nn$ and all $x\in B_{2\delta}(\xk)\setminus B_{\rhok}(\xk)$. Note that $1\not\in I_1\cup...\cup I_q$. We apply Proposition \ref{prop:funda} with $\uk$, $\rhok$ and $\rk$. Similarly to what was done in Step \ref{sec:conc}.2.1, we get, using our induction hypothesis, that
 
$$\lim_{k\to +\infty}\int_{B_{R r_{q+1,k}}(\xk)\setminus \ovb_{R_1 r_{q,k}}(\xk)} \Vk e^{4\uk(x)}\, dx=16\pi^2\,\hbox{Card} (I_{q+1}).$$
The claim then follows from this last equality and \eqref{lim:step:conc}.

\medskip\noindent{\it Step \ref{sec:conc}.2.3:} With Step \ref{sec:conc}.2.2, we get that
\bequa\label{eq:1}
\lim_{k\to +\infty}\int_{B_{R r_{q_0,k}}(\xk)} \Vk e^{4\uk(x)}\, dx=16\pi^2 N,
\eequa
for all $R>0$ large enough. Similarly to what what done in Step \ref{sec:conc}.2.2, there exists $R_0>0$ such that for all $x\in B_{4\delta}(\xk)\setminus \ovb_{R_0 r_{q_0,k}}(\xk)$, we have that
$$|x-\xk| e^{\uk(x)}\leq C\hbox{ and }|x-\xk|^2 |\Delta\uk(x)|\leq C.$$
We apply Proposition \ref{prop:funda} with $\rk=\frac{\delta}{2}$ and $\rhok=R_0 r_{q_0,k}$. We then get that
$$\lim_{k\to +\infty}\int_{B_{\frac{\delta}{2}}(\xk)\setminus B_{R_0 r_{q_0,k}}(\xk)}\Vk e^{4\uk(x)}\, dx=0.$$
This limit, (\ref{eq:1}) and Proposition \ref{estimate:pointwise1} yield
$$\lim_{k\to +\infty}\int_{B_{\delta}(x_0)}\Vk e^{4\uk(x)}\, dx=16\pi^2N.$$
This proves the quantification with $N$ points. We have then proved that ${\bf (H}_{\bf N}{\bf )}$ holds. 

\medskip\noindent In particular, we have proved by induction that ${\bf (H}_{\bf N}{\bf )}$ holds for all $N$.

\medskip\noindent{\bf Step \ref{sec:conc}.3 (Proof of Theorem \ref{th:intro}):} We are now in position to prove the Theorem. We let $\uk$ as in the statement of the Theorem. It follows from Propositions \ref{estimate:pointwise1} that the hypothesis of ${\bf (H}_{\bf N}{\bf )}$ hold in the neighborhood of each of the points of $S_0$. As a consequence, we apply locally ${\bf (H}_{\bf N}{\bf )}$. It then follows that 
$$\Vk e^{4\uk}\rightharpoonup \sum_{i=1}^{N}16\pi^2\alpha_i\delta_{x_i}$$
when $k\to +\infty$. And the proof of Theorem \ref{th:intro} is complete.

\end{document}